\documentclass[11pt]{article}
\usepackage{latexsym}

\font\fancy=eusm10 at 12pt

\newcommand{\ffD}{\mbox{\fancy D}}

\newcommand{\K}[1]{\,${#1}$\,}

\newcommand{\bq}{\begin{equation}}
\newcommand{\eq}{\end{equation}}

\newcommand{\bqr}{\begin{eqnarray}}
\newcommand{\eqr}{\end{eqnarray}}

\newcommand{\bqrx}{\begin{eqnarray*}}
\newcommand{\eqrx}{\end{eqnarray*}}

\newcommand{\br}{\begin{array}}
\newcommand{\er}{\end{array}}

\newcommand{\lcm}{\mbox{lcm}}

\newcommand{\haf}{\scriptsize{\frac{1}{2}}}

\newcommand{\sfrac}[2]{ \mbox{$\frac{#1}{#2}$} }

\newcommand{\voo}{\vspace*{5pt}}

\newcommand{\foo}{\hspace*{4mm}}

\newcommand{\blsk}{\baselineskip}

\begin{document}

\pagestyle{plain}
\pagenumbering{arabic}

\setlength{\parindent}{20pt}
\setlength{\footskip}{.4in}



\setlength{\blsk}{.21in}%
\[\]
\begin{center}
A COMPUTATIONAL APPROACH TO FACTORING LARGE INTEGERS 
\end{center}
\mbox{}
\begin{center}
NELSON PETULANTE\voo\\
Department of Mathematics\\
Bowie State University\\
Bowie, MD 20715\\
USA\voo\voo\voo
\end{center}
\setlength{\blsk}{.27in}%
\[\]
\begin{abstract}
To factor an integer $N$, given that it is equal to the product of two primes $p$ and $q$, it suffices to find an integer ${d<\sfrac{1}{2}{N}}$ satisfying the test ``${C(2\sqrt{Nd})^{2}-4Nd}$ is a perfect square", where $C$ denotes the integer ceiling function. In this approach, the factorization problem equates to the problem of designing an optimal data base $\mbox{\fancy D}$ of values $d$ to be tested.\\
\end{abstract}
1991 AMS SUBJECT CLASSIFICATION CODES. \,\,11D09, 14G05.\\
\mbox{}\\
\setcounter{section}{1}
\noindent
1.{\foo}{THE BASIC TEST}

If \K{x} is a positive real number, let \K{C(x)} denote the smallest integer greater than or equal to \K{x} (the integer ceiling function). Let the function \K{f(x)} be defined by the formula \K{f(x)=C(2\sqrt{x})^{2}-4x} and let \K{f^{2}(x)} denote the composition \K{f(f(x))}.

{PROPOSITION 1.1.} Let \K{n} be a positive integer. Then \K{f(n)=0} if and only if \K{n} is a perfect square (\K{\sqrt{n}} is a whole number).

{PROOF.} If \K{f(n)=0}, then \K{C(2\sqrt{n})^{2}=4n}. Thus, \K{C(2\sqrt{n})=2\sqrt{n}}is an {\em even} integer. Therefore, \K{\sqrt{n}} is an integer. Conversely, if \K{\sqrt{n}} is an integer, then so is \K{2\sqrt{n}}. Thus, \K{C(2\sqrt{n})=2\sqrt{n}}, which means that \K{C(2\sqrt{n})^2=4n} and \K{f(n)=0}. 

{PROPOSITION 1.2.} Let \K{n} be a positive integer. Then \K{f^{2}(n)=0} if and only if there exist positive integers \K{u} and \K{v} such that \K{n=uv} and \K{|\sqrt{u}-\sqrt{v}|\leq 1}.

{PROOF.} Suppose that \K{f^{2}(n)=0}. By Proposition 1.1, \K{f(n)=t^{2}} for some positive integer \K{t}. Thus, \K{4n=C(2\sqrt{n})^{2}-t^{2}}. Note that \K{C(2\sqrt{n})} and \K{t} must be both even or both odd. Therefore, if we let \K{u=\haf(C(2\sqrt{n})+t)} and \K{v=\haf(C(2\sqrt{n})-t)}, then both \K{u} and \K{v} are integers and \K{n=uv}. To show that \K{|\sqrt{u}-\sqrt{v}|\leq 1}, observe that \K{u+v=C(2\sqrt{n})}, so that \K{u+v-2\sqrt{uv}=C(2\sqrt{n})-2\sqrt{n}\leq 1}. That is, \K{|\sqrt{u}-\sqrt{v}|^{2}\leq 1}, which implies \K{|\sqrt{u}-\sqrt{v}|\leq 1}. 
\ \ \ \ Conversely, suppose \K{n=uv} where \K{|\sqrt{u}-\sqrt{v}|\leq 1}. By Proposition 1.1, to show that \K{f^{2}(n)=0}, it suffices to show that \K{f(n)=t^{2}} for some integer \K{t}. If \K{u=v}, then \K{n} is a perfect square and there is nothing to prove. So we may assume that \K{u\neq v}. Let \K{t=|u-v|}. Then \K{t^{2}=(u+v)^{2}-4uv}. Since \K{|\sqrt{u}-\sqrt{v}|\leq 1}, it follows that \K{u-2\sqrt{uv}+v\leq 1}. Therefore, since \K{u\neq v}, we have \K{2\sqrt{uv}\leq u+v\leq 2\sqrt{uv}+1}, so that \K{u+v=C(2\sqrt{n})} and \K{t^{2}=C(2\sqrt{n})^{2}-4n=f(n)}. Thus, \K{f^{2}(n)=0}.

{COROLLARY 1.3.} Suppose \K{N} is the product of two primes \K{p} and \K{q} where \K{|\sqrt{p}-\sqrt{q}|\leq 1}. Then the prime factors can be recovered explicitly in terms of \K{N} by way of the formulas: \K{p=\haf(C(2\sqrt{N})+t)} and \K{q=\haf(C(2\sqrt{N})-t)}, where \K{t=\sqrt{C(2\sqrt{N})^2-4N}}.

{EXAMPLE 1.4.} Take \K{N=176039}. Then \K{t=\sqrt{C(2\sqrt{N})^2-4N}=38}, \K{p=\haf(C(2\sqrt{N})+t)=439}
and \K{q=\haf(C(2\sqrt{N})-t)=401}. Thus, the factorization \K{N=(439)(401)} follows instantly from the fact that \K{\sqrt{439}-\sqrt{401}\leq 1}.

{PROPOSITION 1.5.} To factor an integer \K{N}, given that it is equal to the product of two primes \K{p} and \K{q}, it suffices to find an integer \K{d<\sfrac{1}{2}N} satisfying the test \K{f^{2}(Nd)=0}. Then \K{Nd=uv}, where \K{u=\haf(C(2\sqrt{Nd})+t)}, \K{v=\haf(C(2\sqrt{Nd})-t)} and \K{t=\sqrt{C(2\sqrt{Nd})^2-4Nd}}. The prime factors \K{p} and \K{q} can be recovered separately as factors of \K{u} and \K{v} through the formulas \K{p=\gcd(N,u)} and \K{q=\gcd(N,v)}.

{PROOF.} As in Proposition 1.2, with \K{n=Nd}, we have \K{f^{2}(Nd)=0}. Thus, \K{Nd} factors as \K{Nd=uv} where \K{u=\haf(C(2\sqrt{Nd})+t)}, \K{v=\haf(C(2\sqrt{Nd})-t)} and \K{t=\sqrt{C(2\sqrt{Nd})^{2}-4Nd}}. It remains to show that \K{p} and \K{q} are factors of \K{u} and \K{v} separately. A rough estimate is enough to prove this.  Since, by Proposition 1.2, \K{\sqrt{u}-\sqrt{v}=\delta \leq 1}, we see $v\geq u-2\sqrt{u}$, so that, at least whenever $u\geq 16$, we have \K{Nd=uv\geq u^{2}-2u^{\frac{3}{2}}\geq\haf{u^{2}}}. If \K{p} and \K{q} both were factors of \K{u}, we would have \K{u=apq=aN} for some integer \K{a} which would imply that \K{Nd\geq\haf a^{2}N^{2}} or \K{d\geq \haf a^{2}N}, contradicting the assumption that \K{d<\sfrac{1}{2}N}. A similar contradiction occurs if we suppose that \K{p} and \K{q} both divide \K{v}.

{EXAMPLE 1.6.} Take \K{N=1110757} and \K{d=170}. Then \K{f^{2}(Nd)=0}. Employing the formulas in Proposition 1.5, we get \K{t=23}, \K{u=13753} and \K{v=13730}. Thus \K{N=pq}, where \K{p=\gcd(N,u)=1373}, \K{q=\gcd(N,v)=809}.
 
{PROPOSITION 1.7.} Suppose \K{N=pq} where \K{p} and \K{q} are distinct primes. An integer \K{d<\sfrac{1}{2}{N}} satisfies the test \K{f^{2}(Nd)=0} if and only if \K{d} factors as \K{d=xy} where \K{|\sqrt{py}-\sqrt{qx}|\leq 1}. 

{PROOF.} Suppose \K{d} satisfies the test \K{f^{2}(Nd)=0}. Then, as in Proposition 1.2, \K{Nd=uv} where \K{|\sqrt{u}-\sqrt{v}|\leq 1}. By Proposition 1.5, we may assume that \K{u} is a multiple of \K{p}, say \K{u=py}, and \K{v} is a multiple of \K{q}, say \K{v=qx}. Thus, \K{d=xy} and \K{|\sqrt{py}-\sqrt{qx}|\leq 1}. Conversely, suppose \K{d=xy} is found to satisfy the inequality \K{|\sqrt{py}-\sqrt{qx}|\leq 1}. Then, by Proposition 1.2, \K{Nd=pyqx} satisfies the test \K{f^{2}(Nd)=0}.

\setcounter{section}{2}
\noindent
2.{\foo}A FACTORIZATION STRATEGY

Given \K{N=pq} where the distinct prime factors \K{p} and \K{q} are unknown, our objective is to find an integer \K{d} to satisfy the test \K{f^{2}(Nd)=0}. By Proposition 1.7, the test \K{f^{2}(Nd)=0} is successful if any one of the factorizations of \K{d} as \K{d=xy} satisfies the inequality \K{|\sqrt{py}-\sqrt{qx}|\leq 1}. Thus, for a test value \K{d}, the likelihood of success increases as \K{\tau(d)}, the number of divisors of \K{d}, increases. To formalize this observation, we introduce a new function defined on finite sets of integers called the ``yield function".  

{DEFINITION 2.1.} Let \K{d} be a positive integer. The {\em yield} of \K{d}, denoted \K{Y(d)}, is the number of distinct fractions \K{0<\frac{x}{y}<1} in lowest terms such that \K{xyz^2=d} for some integer \K{z}. If \K{S=\{d_{1}, d_{2},\, .\,.\,.\, , d_{k}\}} is a set of test values, then the yield of \K{S}, denoted \K{Y(S)}, is the number of distinct fractions \K{0<\frac{x}{y}<1} in lowest terms such that \K{xyz^{2}\in S} for some integer \K{z}.

{EXAMPLE 2.2.} Let \K{d=12}. The set of distinct fractions \K{0<\frac{x}{y}<1} such that \K{xyz^{2}=12} is \K{\{\frac{1}{12}, \frac{1}{3}, \frac{3}{4}\}}. Note that the fraction \K{\frac{1}{3}} corresponds to the factorization \K{12=(2)(6)=(1)(3)(2^{2})}. Thus \K{Y(12)=3}. 

{EXAMPLE 2.3.} Let \K{S=\{5,12,20\}}. The set of distinct fractions \K{0<\frac{x}{y}<1} such that \K{xyz^{2}\in S} is \K{\{\frac{1}{20},\frac{1}{12},\frac{1}{5},\frac{1}{3},\frac{3}{4},\frac{4}{5}\}}. Thus \K{Y(S)=6}. Note that the factorization \K{5=(1)(5)} contributes nothing to the yield of \K{S} in view of the factorization \K{20=(2)(10)}. In fact, the yield of \K{S} is the same as the yield of the subset \K{S'=\{12,20\}}. 
 
Now, let \K{\ffD} denote a finite data base of test values, say \K{\ffD=[d_{1}, d_{2},\, .\,.\,.\, , d_{m}]}, structured as a list of integers in ascending order. By definition, the {\em cost} of factoring the integer \K{N} relative to \K{\ffD} is the number of values of \K{d\in \ffD} which need to be tested before a successful value (satisfying \K{f^{2}(Nd)=0}) is found. Obviously, our main objective is to construct a data base which minimizes the cost of factoring any given \K{N} of the form \K{N=pq}. Intuitively, at least, it appears fairly evident that some data bases will be more effective than others. Qualitatively speaking, the cost of factoring \K{N} relative to \K{\ffD} should decrease as the yield \K{Y(\ffD)} increases. 

At this stage, we will attempt to supply at least a rough estimate of the effectiveness of a data base. Suppose \K{N=pq}, where \K{p<q}. In \K{\ffD} we want to find \K{d=xy} such that \K{\sqrt{py}-\sqrt{qx}\leq 1}. This is equivalent to 

\bq \sqrt{\frac{p}{q}}-\sqrt{\frac{x}{y}}\leq \frac{1}{\sqrt{qy}}. \eq 
This inequality is satisfied if the set of fractions \K{\frac{x}{y}} in the interval \K{[0,1]} is so numerous that, for at least one of them, \K{\sqrt{\frac{x}{y}}} comes within a distance of \K{\frac{1}{\sqrt{qy}}} of the fixed quantity \K{\sqrt{\frac{p}{q}}}. This will have a high probability of happening if 

\bq Y(\ffD)\geq \sqrt{qy}. \eq

At this point, we need to formulate a reasonable estimate for \K{\sqrt{qy}}. However, without any specific knowledge of the structure of \K{\ffD} this is difficult to do. We turn then to a discussion of some specific choices of data base \K{\ffD}.

\setcounter{section}{2}
\noindent
3.{\foo} SOME SPECIAL DATA BASES

The simplest data base is a list of consecutive integers starting at 1. Let \K{\ffD_{0}(m)=[1,2,3,\, .\,.\,.\, , m]}. To refine the inequality (2), note that \K{y\leq m}. However, the median divisor of a typical \K{d\in\ffD} is \K{\sqrt{d}}, the maximum of which is \K{\sqrt{m}}. It follows that \K{\sqrt{m}} is a good candidate to represent \K{y} in the inequality \K{Y(\ffD)\geq \sqrt{qy}}. Thus, for the data base  \K{\ffD_{0}(m)}, (2) becomes 

\bq Y(\ffD_{0}(m))\geq \sqrt{q}\sqrt[4]{m}. \eq
The explicit dependence on \K{q} can be removed by setting \K{R=\frac{q}{p}>1}, so that \K{\sqrt{q}=\sqrt[4]{R}\sqrt[4]{N}}. Then (3) becomes

\bq Y(\ffD_{0}(m))\geq \sqrt[4]{mNR}. \eq
A sharp lower bound estimate of \K{Y(\ffD_{0}(m))} is given by \K{m} itself. To see this, note that \K{\sum_{d\in\ffD_{0}(m)} Y(d)} is bounded above by \K{\sum_{1\leq k \leq m} \tau(k) = O(n\ln(n))} \cite{Rose}. Thus, the data base \K{\ffD_{0}(m)} has a good chance of factoring \K{N} if \K{m\geq \sqrt[4]{mNR}}. Equivalently, \K{\ffD_{0}(m)} has a good chance of factoring \K{N} if \K{N\leq \frac{m^{3}}{R}}. Often, a successful value of \K{d\in \ffD_{0}(m)} is found well before the entire data base is exhausted. This is borne out by extensive numerical experiments using MAPLE.  

A more interesting type of data base consists of the set of divisors of a given integer \K{B}. Let \K{\ffD_{1}(B)=[d_{1}, d_{2},\, .\,.\,.\, , d_{m}]}, where the \K{d_{j}} are the divisors of \K{B} arranged in increasing order. Thus, \K{m=\tau(B)} and the largest element in \K{\ffD_{1}(B)} is \K{d_{m}=B}.
Let \K{B=p_{1}^{r_{1}} p_{2}^{r_{2}} p_{3}^{r_{3}} \, ...\, p_{k}^{r_{k}}} (prime power factorization), then
\K{m=(r_{1}+1)(r_{2}+1)\,...\,(r_{k}+1)} and
\K{Y_{1}(B)= Y(\ffD_{1}(B))=(2r_{1}+1)(2r_{2}+1)\,...\,(2r_{k}+1)}.
Thus, a good likelihood exists of factoring \K{N} provided that \K{Y_{1}(B)\geq \sqrt[4](N)\sqrt[4](R)\sqrt[8](B)}.

Some interesting choices for \K{B} (evidenced by extensive numerical experiments using MAPLE):\\
\hspace*{20pt} \K{B=n!}\\
\hspace*{20pt} \K{B=(2)(3)(5)\,...\,(p_{k})} (product of first \K{k} primes).\\
\hspace*{20pt} \K{B=\lcm(1,2,3,\,...,\,m)} (\K{\lcm} of first \K{m} integers).\\

{\bf Acknowledgment:} The author gratefully acknowledges the input of Steve Huntsman who found several significant errors and deficiencies in the previous version of this article. Any remaining errors or deficiencies are due solely to the author's negligence.

\end{document}